\documentclass[11pt,reqno]{amsart}
\usepackage{amssymb,verbatim,enumerate,ifthen}
\usepackage[mathscr]{eucal}
\usepackage[T1]{fontenc}
\def\N{\mathbb{N}}
\def\R{\mathbb{R}}
\def\Q{\mathbb{Q}}
\def\Z{\mathbb{Z}}

\def\P{\mathscr{P}}

\def\<{\langle}
\def\>{\rangle}

\def\sign{\mathop{\mbox{\rm sign}}}
\def\ctg{\mathop{\mbox{\rm ctg}}\nolimits}

\newtheorem{theorem}{Theorem}
\newtheorem*{theorem*}{Theorem}
\def\Thm#1#2{\ifthenelse{\equal{#1}{*}}{\begin{theorem*}#2\end{theorem*}}
             {\begin{theorem}\label{T#1}#2\end{theorem}}}
\newtheorem{Atheorem}{Theorem}

\def\thm#1{Theorem~\ref{T#1}}

\newtheorem{proposition}[theorem]{Proposition}
\newtheorem*{proposition*}{Proposition}
\def\Prp#1#2{\ifthenelse{\equal{#1}{*}}{\begin{proposition*}#2\end{proposition*}}
             {\begin{proposition}\label{P#1}#2\end{proposition}}}

\newtheorem{corollary}[theorem]{Corollary}
\newtheorem*{corollary*}{Corollary}
\def\Cor#1#2{\ifthenelse{\equal{#1}{*}}{\begin{corollary*}#2\end{corollary*}}
             {\begin{corollary}\label{C#1}#2\end{corollary}}}
\newtheorem{Acorollary}[Atheorem]{Corollary}

\def\cor#1{Corollary~\ref{C#1}}

\newtheorem{lemma}[theorem]{Lemma}
\newtheorem*{lemma*}{Lemma}
\def\Lem#1#2{\ifthenelse{\equal{#1}{*}}{\begin{lemma*}#2\end{lemma*}}
             {\begin{lemma}\label{L#1}#2\end{lemma}}}
\newtheorem{Alemma}[Atheorem]{Lemma}

\def\lem#1{Lemma~\ref{L#1}}

\theoremstyle{definition}
\newtheorem{remark}[theorem]{Remark}
\newtheorem*{remark*}{Remark}
\def\Rem#1#2{\ifthenelse{\equal{#1}{*}}{\begin{remark*}\rm #2\end{remark*}}
             {\begin{remark}\label{R#1}\rm #2\end{remark}}}

\newtheorem{example}[theorem]{Example}
\newtheorem*{example*}{Example}
\def\Exa#1#2{\ifthenelse{\equal{#1}{*}}{\begin{example*}\rm #2\end{example*}}
             {\begin{example}\label{Ex#1}\rm #2\end{example}}}

\def\eq#1{{\rm(\ref{E#1})}}
\def\Eq#1#2{\ifthenelse{\equal{#1}{*}}
  {\begin{equation*}\begin{aligned}#2\end{aligned}\end{equation*}}
  {\begin{equation}\begin{aligned}\label{E#1}#2\end{aligned}\end{equation}}}

\begin{document}
\begin{flushright}
%\textit{Submitted to: }
\end{flushright}
\vspace{5mm}

%\date{\today}

\title[A new characterization of convexity]
{A new characterization of convexity with respect to Chebyshev systems}

\author[Zs. P\'ales]{Zsolt P\'ales}
\author[\'E. Sz\'ekelyn\'e Rad\'acsi]{\'Eva Sz\'ekelyn\'e Rad\'acsi}
\address{Institute of Mathematics, University of Debrecen,
H-4032 Debrecen, Egyetem t\'er 1, Hungary}
\email{\{pales,radacsi.eva\}@science.unideb.hu}


\thanks{The research was supported by the Hungarian Scientific Research Fund (OTKA) Grant K-111651 and by the 
EFOP-3.6.1-16-2016-00022 project. This project is co-financed by the European Union and the European Social Fund.}

%\dedicatory{Dedicated to the 70th birthday of Professor Roman Ger}

\begin{abstract}
The notion of $n$th order convexity in the sense of Hopf and Popoviciu is defined via the 
nonnegativity of the $(n+1)$st order divided differences of a given real-valued function. 
In view of the well-known recursive formula for divided differences, the nonnegativity of 
$(n+1)$st order divided differences is equivalent to the $(n-k-1)$st order convexity of the
$k$th order divided differences which provides a characterization of $n$th order convexity.

The aim of this paper is to apply the notion of higher-order divided differences in the 
context of convexity with respect to Chebyshev systems introduced by Karlin in 1968. Using a 
determinant identity of Sylvester, we then establish a formula for the generalized divided 
differences which enables us to obtain a new characterization of convexity 
with respect to Chebyshev systems. Our result generalizes that of Wąsowicz which was obtained in 2006. 
As an application, we derive a necessary condition for functions which can 
be written as the difference of two functions convex with respect to a given Chebyshev system. 
\end{abstract}

\maketitle

\section{Introduction}

Denote by $\N$, $\Z$, $\Q$, $\R$ the sets of natural, integer, rational, and real numbers, respectively.
Given a set $H\subseteq\R$, the set of positive elements of $H$ is denoted by $H_+$. Thus, for instance,
$\N=\Z_+$.

For a set $H\subseteq\R$, denote the simplex of strictly increasingly ordered $n$-tuples of 
$H^n$ by $\sigma_n(H)$, i.e.,
\Eq{*}{
  \sigma_n(H):=\{(x_1,\dots,x_n)\in H^n\mid x_1<\dots<x_n\}.
}
The set of those elements of $H^n$ that have pairwise distinct 
coordinates will be denoted by $\tau_n(H)$, i.e.,
\Eq{*}{
  \tau_n(H):=\{ (x_1,\dots,x_n)\in H^n\mid x_i\neq x_j\, (i\neq j)\}.
}
Obviously, $\sigma_n(H)\neq\emptyset$ and $\tau_n(H)\neq\emptyset$ if and only if the cardinality $|H|$ of 
$H$ is at least $n$, furthermore, we have that $\sigma_n(H)\subseteq \tau_n(H) \subseteq H^n$.

Provided that $|H|\geq n$, for a vector valued function $ \omega=(\omega_1, \dots, \omega_n): H \to \R^n $,
the functional operator $\Phi_\omega:H^n\to\R$ is defined by
\Eq{*}{
    \Phi_\omega(x_1,\dots,x_n):=
   \left|\begin{array}{ccc}
     \omega_1(x_1) &  \dots & \omega_1(x_n)\\
     \vdots & \ddots &  \vdots \\
     \omega_n(x_1) &  \dots & \omega_n(x_n)
   \end{array}\right|
   \qquad\big((x_1,\dots,x_n)\in H^n\big).
}
We say that $ \omega $ is an \emph{$n$-dimensional positive (resp.\ negative) Chebyshev system over $H$} if
$\Phi_\omega$ is positive (resp.\ negative) over $\sigma_n(H)$, respectively.

The following systems are the most important particular cases for positive Chebyshev systems. For more
important examples we refer to the books by Karlin \cite{Kar68} and Karlin--Studden \cite{KarStu66}.
\begin{enumerate}[(i)]
\item The function $\omega:\R\to\R^n$ given by $\omega(x):=(1, x, \dots, x^{n-1}) $ is an
$n$-dimensional positive Chebyshev system on $ \R $. This system is called the
\emph{standard, or polynomial $n$-dimensional Chebyshev system}.
\item The function $ \omega(x):=\big( 1,\cos(x), \sin(x), \dots, \cos(nx), \sin(nx)\big) $ is a
$(2n+1)$-dimensional positive Chebyshev system on any open interval $I$ whose length is less than or equal to 
$ 2\pi $.
\item The function $ \omega(x):=\big( \cos(x), \sin(x), \dots, \cos(nx), \sin(nx)\big) $
is a $(2n)$-dimensional positive Chebyshev system on any open interval $I$ whose length is less than or equal
to $ \pi $. 
\item For the function $ \omega(x):=(1,x^2 )$, we get that $\Phi_\omega(x_1,x_2)=(x_1+x_2)(x_2-x_1)$.
Therefore $\omega$ is a $2$-dimensional positive Chebyshev system on $ \R_+ $, but it is not a Chebyshev
system on $ \R $ (observe that $ \Phi_{\omega}(-1,1)=0$).
\end{enumerate}

Given a positive Chebyshev system $\omega:H\to\R^n$, a function $f:H\to\R$ is called
\emph{$\omega$-convex} (i.e., \emph{convex with respect to the Chebyshev system $\omega$}) if,
$\Phi_{(\omega,f)}$ is nonnegative over $\sigma_{n+1}(H)$. A function $f:H\to\R$ is strictly $\omega$-convex
if a function $ (\omega_1,\dots,\omega_n,f)$ is an $(n+1)$-dimensional positive Chebyshev system over $H$.

For $k\geq0$, define the $k$th power function $p_k:\R\to\R$ by $p_k(x):=x^k$. As we have seen it before,
$(p_0,\dots,p_{n-1})$ is an $n$-dimensional Chebyshev system. The notion of convexity with respect to this
system, called polynomial convexity, was introduced by Hopf \cite{Hop26} and by Popoviciu
\cite{Pop44}. The particular case, when $\omega=(p_0,p_1)$, simplifies to the notion of standard convexity,
moreover, for $(x,y,z)\in\sigma_3(H)$, the inequality $\Phi_{(p_0,p_1,f)}(x,y,z)\geq0$ is equivalent to
\Eq{*}{
 f(y)\leq \frac{z-y}{z-x}f(x)+\frac{y-x}{z-x}f(z).
}
It is easy to verify that this inequality holds if and only if, for all $p\in H$, the mapping
\Eq{*}{
  x\mapsto \frac{f(x)-f(p)}{x-p}\qquad(x\in H\setminus\{p\})
}
is nondecreasing. More generally, in view of the well-known recursive formula for divided differences, the 
nonnegativity of $(n+1)$st order divided differences is equivalent to the monotonicity of the
$n$th order divided differences which provides a characterization of $n$th order convexity.

The aim of this paper is to characterize higher-order convexity with respect to Chebyshev systems (cf. 
\cite{BesPal04}) by applying the notion of related generalized higher-order divided differences introduced by 
Karlin in \cite{Kar68} and rediscovered in \cite{PalRad16}. Using a determinant identity of Sylvester, we 
then establish a formula for the generalized divided differences of higher order and obtain new characterizations of 
convexity with respect to Chebyshev systems. Our result generalizes that of Wąsowicz which was obtained in 
\cite{Was06b}. As an application, we introduce the notion of $\omega$-variation and we derive a necessary condition for 
functions that can be written as the difference of two $\omega$-convex functions.

\section{Characterizations of convexity with respect to Chebyshev systems}

In the sequel, we will need the following classical formula which is termed Sylvester's Determinant Identity 
in the literature \cite{AkrAkrMal96}, \cite{Gan59}.

\Thm{*}{Let $n\in\N$ and $A:\{1,\dots,n\}\times\{1,\dots,n\}\to\R$ be an 
$n\times n$ matrix. For $1\leq k \leq n-1$ define
the $(n-k)\times(n-k)$ matrix $B_k:\{k+1,\dots,n\}\times\{k+1,\dots,n\}\to\R$ by
\Eq{*}{
  B_k(i,j):=\det\big(A|_{\{1,\dots,k,i\}\times\{1,\dots,k,j\}}\big) 
\qquad(i,j\in\{k+1,\dots,n\}).
}
Then the following identity holds
\Eq{*}{
   \det(B_k)=\Big(\det\big(A|_{\{1,\dots,k\}\times\{1,\dots,k\}}\big) \Big)^{n-k-1} \det(A).
}}

To formulate our main results below, we introduce the notions of divided differences with respect to 
Chebyshev systems. Let $n\in\N$, $H\subseteq\R$ with $|H|\geq n$ and let
$\omega=(\omega_1,\dots, \omega_n):H\to\R^n$ be an $n$-dimensional positive Chebyshev system over $H$. 
In the sequel, for the sake of convenience 
and brevity, given $k\in\{1,\dots,n\}$, we shall denote by $\omega_{\<k\>}$ the $k$-tuple $(\omega_1,\dots, 
\omega_k)$. Thus, in particular, we have that 
\Eq{*}{
\omega_{\<1\>}=\omega_1,\qquad \omega_{\<2\>}=(\omega_1,\omega_2),
\qquad \dots, \qquad \omega_{\<n\>}=\omega.
}
For a function $f:H\to\R$ and $k\in\{1,\dots,n\}$, the generalized \emph{$(k-1)$-st order $\omega$-divided 
difference of $f$} (cf.\ \cite{Kar68}) is defined by
    \Eq{*}{
        \big[x_1,\dots,x_{k};f\big]_{\omega_{\<k\>}}
        :=\frac{\Phi_{(\omega_{\<k-1\>},f)}(x_1,\dots,x_{k})}{\Phi_{\omega_{\<k\>}}(x_1,\dots,x_{k})}
             \qquad \big( (x_1,\dots,x_{k}) \in \tau_{k}(H) \big).
    }
provided that $\omega_{\<k\>}$ is a $k$-dimensional Chebyshev system. Clearly, if $\omega=(\omega_1,\dots, 
\omega_n)=(p_0,\dots,p_{n-1})$, then, $\big[x_1,\dots,x_k;f\big]_{\omega_{\<k\>}}$ is equal to the standard
$(k-1)$-st order divided difference $\big[x_1,\dots,x_k;f\big]$.

\Thm{1}{Let $n,k\in\N$, $k<n$, $|H| \geq n$ and let $x_1<\dots<x_k$ be arbitrary elements of $H$. Let 
$\omega:=(\omega_1,\dots,\omega_n)$ be an $n$-dimensional positive Chebyshev system over $H$ such that 
$\omega_{\<k\>}$ and $\omega_{\<k+1\>}$ are $k$ and $(k+1)$-dimensional positive Chebyshev system over $H$, 
respectively. Then the following system of functions
\Eq{fs}{
  x \mapsto  \big[x_1,\dots,x_k,x;\omega_j \big]_{\omega_{\<k+1\>}}
  \qquad\qquad\qquad (k+1 \leq j \leq n)
}
is an $(n-k)$-dimensional positive Chebyshev system over $H\setminus\{x_1,\dots,x_k \}$.}

\begin{proof}
Let $x_{k+1}<\dots<x_n$ be arbitrary elements of $H\setminus \{x_1,\dots,x_k \}$.
Applying Sylvester's determinant identity for the matrix $A$ defined by $A(i,j):=\omega_i(x_j)$, we get 
\Eq{*}{
  \Phi_{\omega}(x_1,\dots,x_n) &\cdot \big( \Phi_{\omega_{\<k\>}}(x_1,\dots,x_{k})\big)^{n-k-1} \\
  &= \left|\begin{array}{ccc}
    \Phi_{(\omega_{\<k\>},\omega_{k+1})}(x_1,\dots,x_{k},x_{k+1}) & \dots &
	\Phi_{(\omega_{\<k\>},\omega_{k+1})}(x_1,\dots,x_{k},x_{n}) \\ 
    \vdots & \ddots & \vdots \\
    \Phi_{(\omega_{\<k\>},\omega_{n})}(x_1,\dots,x_{k},x_{k+1}) & \dots &
	\Phi_{(\omega_{\<k\>},\omega_{n})}(x_1,\dots,x_{k},x_{n})
 \end{array}\right|.
}
Then dividing the $j$-th column $(j\in\{1,\dots,n-k\})$ of the determinant on the right hand side of the above 
identity by $ \Phi_{\omega_{\<k+1\>}}(x_1,\dots,x_{k},x_j) $, we arrive at the following equality:
\Eq{id}{
    &\frac{\Phi_{\omega}(x_1,\dots,x_n) \cdot 
    \big( \Phi_{\omega_{\<k\>}}(x_1,\dots,x_{k}) \big)^{n-k-1}}{
    \prod_{j=k+1}^{n} \Phi_{\omega_{\<k+1\>}}(x_1,\dots,x_{k},x_j)}\\
    &\hspace{3cm}=\left|\begin{array}{ccc}
    \big[x_1,\dots,x_k,x_{k+1};\omega_{k+1} \big]_{\omega_{\<k+1\>}} & \cdots 
	  & \big[x_1,\dots,x_k,x_n;\omega_{k+1} \big]_{\omega_{\<k+1\>}} \\
    \vdots & \ddots & \vdots \\
    \big[x_1,\dots,x_k,x_{k+1};\omega_{n} \big]_{\omega_{\<k+1\>}} & \cdots
	  & \big[x_1,\dots,x_k,x_n;\omega_{n} \big]_{\omega_{\<k+1\>}}
    \end{array}\right|.
}
In order to complete the proof of the theorem, it suffices to show that the left hand side of 
the above identity is positive for elements $x_{k+1}<\dots<x_n$ of $H\setminus \{x_1,\dots,x_k \}$.

Define the indices $\ell_{k+1},\dots,\ell_n$ by
\Eq{*}{
  \ell_j:=\begin{cases}
           \max\big\{i\in\{1,\dots,k\}\mid x_i<x_j\big\} &\mbox{if } x_1<x_j, \\[2mm]
           0 &\mbox{if } x_j<x_1.
          \end{cases}
}
Now, using that $\omega_{\<k+1\>}$ is a positive Chebyshev system, for 
$j\in\{k+1,\dots,n\}$, we are going to show that
\Eq{sg1}{
  \sign \Phi_{\omega_{\<k+1\>}}(x_1,\dots,x_{k},x_j)=(-1)^{k-\ell_j}.
}
If $x_j<x_1$, then $\ell_j=0$ and 
\Eq{*}{
\sign\Phi_{\omega_{\<k+1\>}}(x_1,\dots,x_{k},x_j)
 =(-1)^k\sign\Phi_{\omega_{\<k+1\>}}(x_j,x_1,\dots,x_{k})=(-1)^k=(-1)^{k-\ell_j}.
}
If $x_1<x_j<x_k$, then $x_{\ell_j}<x_j<x_{\ell_j+1}$, hence
\Eq{*}{
\sign\Phi_{\omega_{\<k+1\>}}(x_1,\dots,x_{k},x_j)
=(-1)^{k-\ell_j}\sign\Phi_{\omega_{\<k+1\>}}(x_1,\dots,x_{\ell_j},x_j,x_{\ell_j+1},\dots,x_{k})
=(-1)^{k-\ell_j }.
}
Finally, if $x_k<x_j$, then $\ell_j=k$ and
\Eq{*}{
\sign\Phi_{\omega_{\<k+1\>}}(x_1,\dots,x_{k},x_j)=1=(-1)^{k-\ell_j }.
}
Applying \eq{sg1}, we get
\Eq{*}{
  \sign \prod_{j=k+1}^{n} 
\Phi_{\omega_{\<k+1\>}}(x_1,\dots,x_{k},x_j)=(-1)^{(n-k)k-(\ell_{k+1}+\cdots+\ell_n)}.
}
An analogous computation and the positive Chebyshev property of $\omega$ results that
\Eq{*}{
  \sign \Phi_{\omega}(x_1,\dots,x_n)=(-1)^{(n-k)k-(\ell_{k+1}+\cdots+\ell_n)}.
}
Therefore, the left hand side of \eq{id} is positive since $\omega_{\<k\>}$
is also a positive Chebyshev system.
\end{proof}

\Thm{2}{Let $n,k\in\N$, $k<n$, $|H| \geq n+1$. Let $\omega:=(\omega_1,\dots,\omega_n):H\to\R^n$ be an 
$n$-dimensional positive Chebyshev system over $H$ such that 
$\omega_{\<k\>}$ and $\omega_{\<k+1\>}$ are $k$ and $(k+1)$-dimensional positive Chebyshev system over $H$, 
respectively and let $f:H\to \R$ be a function.
Then, for all $ (x_1,\dots,x_{n+1}) \in \tau_{n+1}(H)$, the following identity is valid
\Eq{fid}{
    &\frac{\Phi_{(\omega,f)}(x_1,\dots,x_{n+1}) \cdot 
    \big( \Phi_{\omega_{\<k\>}}(x_1,\dots,x_{k}) \big)^{n-k}}{
    \prod_{j=k+1}^{n+1} \Phi_{\omega_{\<k+1\>}}(x_1,\dots,x_{k},x_j)}\\
    &\hspace{3cm}=\left|\begin{array}{cccc}
    \big[x_1,\dots,x_k,x_{k+1};\omega_{k+1} \big]_{\omega_{\<k+1\>}} & \cdots 
	  & \big[x_1,\dots,x_k,x_{n+1};\omega_{k+1} \big]_{\omega_{\<k+1\>}} \\
    \vdots & \ddots & \vdots \\
    \big[x_1,\dots,x_k,x_{k+1};\omega_{n} \big]_{\omega_{\<k+1\>}} & \cdots
	  & \big[x_1,\dots,x_k,x_{n+1};\omega_{n} \big]_{\omega_{\<k+1\>}} \\
    \big[x_1,\dots,x_k,x_{k+1};f \big]_{\omega_{\<k+1\>}} & \cdots
	  & \big[x_1,\dots,x_k,x_{n+1};f \big]_{\omega_{\<k+1\>}}
    \end{array}\right|.
}
Furthermore, the following statements are equivalent:
\begin{enumerate}[(i)]
 \item $f$ is $\omega$-convex on $H$; 
 \item For each ordered $k$-tuple $(x_1,\dots,x_k)\in\sigma_k(H)$, the function
$ x \mapsto\big[x_1,\dots,x_k,x;f \big]_{\omega_{\<k+1\>}} $ is convex on $H\setminus\{x_1,\dots,x_k\}$ with 
respect to the $(n-k)$-dimensional Chebyshev system defined by \eq{fs};
 \item There exists $\ell\in\{0,\dots,k\}$ such that, for each ordered $k$-tuple 
$(x_1,\dots,x_k)\in\sigma_k(H)$, the function $ x \mapsto\big[x_1,\dots,x_k,x;f \big]_{\omega_{\<k+1\>}} $ is 
convex with respect to the $(n-k)$-dimensional Chebyshev system defined by 
\eq{fs} on $H\cap\,]-\infty,x_{1}[$ if $\ell=0$, on $H\cap\,]x_\ell,x_{\ell+1}[$ if $0<\ell<k$ and on 
$H\cap\,]x_k,+\infty[$ if $\ell=k$.
\end{enumerate}}

\begin{proof}
The formula in \eq{fid} follows from Sylvester's determinant identity with the 
$(n+1)\times (n+1)$ matrix $A$ defined by
\Eq{*}{
    A(i,j):= \begin{cases}
               \omega_i(x_j) \qquad\mbox{if } (i,j)\in\{1,\dots,n\}\times\{1,\dots,n+1\},\\
               f(x_j)  \qquad\mbox{ if } (i,j)\in\{n+1\}\times\{1,\dots,n+1\}
             \end{cases}
}
in the same way as formula \eq{id} in the proof of \thm{1}.

(i)$\Rightarrow$(ii). Assume that the function $f$ is $\omega$-convex, i.e. $ \Phi_{(\omega,f)}$ is 
nonnegative over $\sigma_{n+1}(H) $. 
Let $ k<n $, let $ x_1,\dots,x_k \in \sigma_k(H) $ and $ x_{k+1},\dots,x_{n+1}\in \sigma_{n+1-k}(H\setminus\{x_1,\dots,x_k\}) $.
With similar idea as in the previous proof we can prove that the left hand side
of \eq{fid} is nonnegative, hence the right hand side of \eq{fid} is also nonnegative, then using the positive 
Chebyshev property of \eq{fs}, we get that a function $ x \mapsto\big[x_1,\dots,x_k,x;f \big]_{\omega_{\<k+1\>}} $ 
is a convex function with respect to the Chebyshev system \eq{fs}.

The implication (ii)$\Rightarrow$(iii) is trivial.

(iii)$\Rightarrow$(i). Assume that (iii) holds for some $\ell\in\{0,\dots,k\}$. To prove that $f$ is 
$\omega$-convex, let $(x_1,\dots,x_{n+1})\in\sigma_{n+1}(H)$. Define 
\Eq{*}{
  x'_i:=\begin{cases}
        x_i &\mbox{ if } 1\leq i \mbox{ and } i\leq\ell,\\[2mm]  
        x_{i+n-k+1} &\mbox{ if } \ell+1\leq i \mbox{ and } i\leq k,\\[2mm] 
        x_{i+\ell-k} &\mbox{ if } k+1\leq i \mbox{ and } i\leq n+1. 
        \end{cases}
}
Now observe that $(x'_1,\dots,x'_{k})\in\sigma_k(H)$ and $(x'_{k+1},\dots,x'_{n+1})\in 
\sigma_{n-k+1}(H_\ell)$, where 
\Eq{*}{
  H_\ell:=\begin{cases}
          H\cap\,]-\infty,x'_{1}[ &\mbox{ if } \ell=0,\\[2mm] 
          H\cap\,]x'_\ell,x'_{\ell+1}[ &\mbox{ if } 0<\ell<k,\\[2mm]
          H\cap\,]x'_{k},\infty[ &\mbox{ if } \ell=k.
          \end{cases}
}
To complete the proof, applying formula \eq{fid} for the $(n+1)$-tuple $(x'_1,\dots,x'_{n+1})$, we obtain:
\Eq{fid'}{
    &\frac{\Phi_{(\omega,f)}(x'_1,\dots,x'_{n+1}) \cdot 
    \big( \Phi_{\omega_{\<k\>}}(x'_1,\dots,x'_{k}) \big)^{n-k}}{
    \prod_{j=k+1}^{n+1} \Phi_{\omega_{\<k+1\>}}(x'_1,\dots,x'_{k},x'_j)}\\
    &\hspace{3cm}=\left|\begin{array}{cccc}
    \big[x'_1,\dots,x'_k,x'_{k+1};\omega_{k+1} \big]_{\omega_{\<k+1\>}} & \cdots 
	  & \big[x'_1,\dots,x'_k,x'_{n+1};\omega_{k+1} \big]_{\omega_{\<k+1\>}} \\
    \vdots & \ddots & \vdots \\
    \big[x'_1,\dots,x'_k,x'_{k+1};\omega_{n} \big]_{\omega_{\<k+1\>}} & \cdots
	  & \big[x'_1,\dots,x'_k,x'_{n+1};\omega_{n} \big]_{\omega_{\<k+1\>}} \\
    \big[x'_1,\dots,x'_k,x'_{k+1};f \big]_{\omega_{\<k+1\>}} & \cdots
	  & \big[x'_1,\dots,x'_k,x'_{n+1};f \big]_{\omega_{\<k+1\>}}
    \end{array}\right|.
}
We have that the right hand side of the above equality is nonnegative, since, by (iii), 
$\big[x'_1,\dots,x'_k,\cdot,f\big]_{\omega_{\<k+1\>}}$ is convex with respect to the $(n-k)$-dimensional 
Chebyshev system defined by \eq{fs} (where the $x_i$s are replaced by $x'_i$) on $H_\ell$. The subsystem 
$\omega_{\<k\>}$ being a positive Chebyshev system, $\Phi_{\omega_{\<k\>}}(x'_1,\dots,x'_{k})>0$. Hence,  
\eq{fid'} implies that
\Eq{a}{
  \frac{\Phi_{(\omega,f)}(x'_1,\dots,x'_{n+1})}{
    \prod_{j=k+1}^{n+1} \Phi_{\omega_{\<k+1\>}}(x'_1,\dots,x'_{k},x'_j)}\geq0.
}
For $j\in\{k+1,\dots,n+1\}$, we have that $x'_j<x'_1$ if $\ell=0$,  $x'_\ell<x'_j<x'_{\ell+1}$ if $0<\ell<k$, 
and $x'_k<x'_j$ if $\ell=k$, therefore
\Eq{*}{
  \sign \Phi_{\omega_{\<k+1\>}}(x'_1,\dots,x'_{k},x'_j)=(-1)^{k-\ell}, 
}
which yields that
\Eq{b}{
  \sign \prod_{j=k+1}^{n+1} \Phi_{\omega_{\<k+1\>}}(x'_1,\dots,x'_{k},x'_j) = (-1)^{(n-k+1)(k-\ell)}.
}
Therefore, using \eq{b} and inequality \eq{a}, after interchanging the appropriate columns of the 
determinant $\Phi_{(\omega,f)}(x_1,\dots,x_{n+1})$, we get that
\Eq{*}{
  \Phi_{(\omega,f)}(x_1,\dots,x_{n+1})=(-1)^{(n-k+1)(k-\ell)}\Phi_{(\omega,f)}(x'_1,\dots,x'_{n+1})\geq0.
}
This completes the proof of the $\omega$-convexity of $f$.
\end{proof}

The following result, which was established by Wąsowicz \cite[Theorem 2]{Was06b}, concerns the particular 
case $k=n-1$ of the previous theorem.

\Cor{1}{Let $n\in\N$, $n\geq 2$, $|H| \geq n+1$ and $\omega:=(\omega_1,\dots,\omega_n):H\to \R^n$ be an 
$n$-dimensional positive Chebyshev system such that $\omega_{\<n-1\>}$ is an $(n-1)$-dimensional positive
Chebyshev system and let $f:H\to \R$ be a function. Then, for all $ (x_1,\dots,x_{n+1}) \in \tau_{n+1}(H)$, 
the following identity is valid
\Eq{fid1}{
    \frac{\Phi_{(\omega,f)}(x_1,\dots,x_{n+1})
    \Phi_{\omega_{\<n-1\>}}(x_1,\dots,x_{n-1})}{
    \Phi_{\omega}(x_1,\dots,x_{n-1},x_n)\Phi_{\omega}(x_1,\dots,x_{n-1},x_{n+1})}
    =\big[x_1,\dots,x_{n-1},x_{n+1};f \big]_{\omega}-\big[x_1,\dots,x_{n-1},x_{n};f \big]_{\omega}.
}
Furthermore, the following statements are equivalent:
\begin{enumerate}[(i)]
 \item $f$ is $\omega$-convex;
 \item For each ordered $(n-1)$-tuple $ (x_1,\dots,x_{n-1})\in\sigma_{n-1}(H) $, 
	  the function $x \mapsto \big[x_1,\dots,x_{n-1},x;f\big]_{\omega}$ is nondecreasing on 
	  $H \setminus \{x_1,\dots,x_{n-1} \}$;
 \item There exists $\ell\in\{0,\dots,n-1\}$ such that, for each ordered $(n-1)$-tuple
	  $(x_1,\dots,x_{n-1})\in\sigma_{n-1}(H) $, 
	  the function $x \mapsto \big[x_1,\dots,x_{n-1},x;f\big]_{\omega}$ is nondecreasing on
	  $H\cap\,]-\infty,x_{1}[$ if $\ell=0$, on $H\cap\,]x_{\ell},x_{\ell+1}[$ if $0<\ell<n-1$ and on 
	  $H\cap\,]x_{n-1},+\infty[$ if $\ell=n-1$.
\end{enumerate}}

\begin{proof} If $k=n-1$ then the $n-k=1$ dimensional Chebyshev system defined by \eq{fs} is the constant 
function $1$ and \eq{fid1} is a particular case of \eq{fid}. The convexity of the function $x \mapsto 
\big[x_1,\dots,x_{n-1},x;f\big]_{\omega}$ with respect to this Chebyshev system is equivalent to its 
nondecreasingness. Thus, \thm{2} directly yields the equivalence of statements (i), (ii), and (iii).
\end{proof}

In what follows, we apply \thm{2} to the $n$-dimensional polynomial system. For this, we shall need the 
following auxiliary statement. Recall that we have defined $p_n:\R\to\R$ by $p_n(x):=x^{n}$. 

\Lem{3}{Let $k\in\N$ and $x_1<\dots<x_{k}$ be arbitrary elements of $H$. Then the following 
equality is valid 
\Eq{pi}{
    \big[x_1,\dots,x_{k};p_n \big] =
    \sum_{\substack{\alpha_1,\dots,\alpha_{k} \geq 0,\\ \alpha_1+\dots+\alpha_{k}=n-k+1}}
    x_1^{\alpha_1} \cdots x_{k}^{\alpha_{k}}
    \qquad\qquad (n\in\N\cup\{0\}).
}}

\begin{proof}
The proof runs by induction on $k$. For $k=1$ the statement trivially holds. Assume that
\eq{pi} is true for $k=m-1 \in \N$, $m \geq 2$. By a well-known property of classical divided differences, we 
have
\Eq{*}{
    \big[x_1,\dots,x_{m};p_n \big]=
    \frac{\big[x_2,\dots,x_{m};p_n \big]-\big[x_1,\dots,x_{m-1};p_n \big]}
	{x_m-x_1}.
}
By the induction hypothesis,
\Eq{*}{
    \big[x_1,\dots,x_{m};p_n \big] &=
    \displaystyle\sum_{\substack{\alpha_2,\dots,\alpha_{m} \geq 0,\\ \alpha_2+\dots+\alpha_{m}=n-m+2}}
    \frac{x_2^{\alpha_2} \cdots x_{m}^{\alpha_{m}}}{x_{m}-x_1}-
    \displaystyle\sum_{\substack{\alpha_1,\dots,\alpha_{m-1} \geq 0,\\ \alpha_1+\dots+\alpha_{m-1}=n-m+2}}
    \frac{x_1^{\alpha_1} \cdots x_{m-1}^{\alpha_{m-1}}}{x_{m}-x_1}\\
    &= \sum_{j=0}^{n-m+2} \Bigg({\frac{x_{m}^j - x_1^j}{x_{m}-x_1}}\!\!\!
    \sum_{\substack{\alpha_2,\dots,\alpha_{m-1} \geq 0,\\ \alpha_2+\dots+\alpha_{m-1}=n-m-j+2}} \!\!\!
    {x_2^{\alpha_2} \cdots x_{m-1}^{\alpha_{m-1}}}\Bigg)%\\
    = \sum_{\substack{\alpha_1,\dots,\alpha_{m} \geq 0,\\ \alpha_1+\dots+\alpha_{m}=n-m+1}}
    {x_1^{\alpha_1} \cdots x_{m}^{\alpha_{m}}}.
}
Thus we obtain \eq{pi} for $k=m$, which completes the proof of the induction.
\end{proof}

\Cor{2}{Let $n,k\in\N$, $k<n$, $|H| \geq n+1$. Let $ f:H\to\R $ be a function. Then the following statements 
are pairwise equivalent.
\begin{enumerate}[(i)]
 \item $f$ is $n$-monotone, i.e., it is convex with respect to the Chebyshev system $(p_0,\dots,p_{n-1})$; 
 \item For each ordered $k$-tuple $ (x_1,\dots,x_{k})\in \sigma_{k}(H) $, the 
function $x\mapsto \big[ x_1,\dots,x_k,x;f \big]$ is $(n-k)$-monotone (i.e., it is convex with respect to the 
Chebyshev system $(p_0,\dots,p_{n-k-1})$) on the set $H\setminus\{ x_1,\dots,x_k \} $;
 \item There exists $\ell\in\{0,\dots,k\}$ such that, for each ordered $k$-tuple 
$(x_1,\dots,x_k)\in\sigma_k(H)$, the function $ x \mapsto\big[x_1,\dots,x_k,x;f \big] $ is $(n-k)$-monotone 
on $H\cap\,]-\infty,x_{1}[$ if $\ell=0$, on $H\cap\,]x_\ell,x_{\ell+1}[$ if $0<\ell<k$ and on 
$H\cap\,]x_k,+\infty[$ if $\ell=k$.
\end{enumerate}
}

\begin{proof} Let $ (x_1,\dots,x_{k})\in \sigma_{k}(H) $ be fixed. Define, for $\ell\geq0$,
\Eq{*}{
  P_\ell:=P_\ell(x_1,\dots,x_k)
  :=\sum_{\substack{\alpha_1,\dots,\alpha_{k} \geq 0,\\ \alpha_1+\dots+\alpha_{k}=\ell}}
    {x_1^{\alpha_1} \cdots x_{k}^{\alpha_{k}}}.
}
Observe that $P_0=1$, $P_1=x_1+\cdots+x_k$, etc. Using \lem{3}, for $j\in\{k,\dots,n-1\}$, we obtain
\Eq{*}{
  \big[x_1,\dots,x_k,x;p_j \big]
  =\sum_{\alpha=0}^{j-k}\Bigg(\sum_{\substack{\alpha_1,\dots,\alpha_{k} \geq 0,\\ 
\alpha_1+\dots+\alpha_{k}=j-k-\alpha}}
    {x_1^{\alpha_1} \cdots x_{k}^{\alpha_{k}}}\Bigg)x^\alpha
  =\sum_{\alpha=0}^{j-k} P_{j-k-\alpha}x^\alpha.
}
Therefore, performing elementary row operations on determinants (subtracting $P_1$ times the first row from 
the second, then subtracting $P_2$ times the first plus $P_1$ times the second row from the third, etc.), for 
the right hand side of \eq{fid}, we obtain the following formula
\Eq{*}{
    &\hspace{-10mm}\left|\begin{array}{lll}
       \big[x_1,\dots,x_k,x_{k+1};p_k \big] & \cdots 
	  & \big[x_1,\dots,x_k,x_{n+1};p_k \big] \\[1mm]
       \qquad\qquad\vdots & \ddots & \qquad\qquad\vdots \\
       \big[x_1,\dots,x_k,x_{k+1};p_{n-1} \big] & \cdots
	  & \big[x_1,\dots,x_k,x_{n+1};p_{n-1} \big] \\[1mm]
       \big[x_1,\dots,x_k,x_{k+1};f \big]& \cdots
	  & \big[x_1,\dots,x_k,x_{n+1};f \big]
    \end{array}\right|\\[2mm]
    &=\left|\begin{array}{ccc}
       1 & \cdots & 1 \\
       x_{k+1}+P_1& \cdots & x_{n+1}+P_1 \\
       \vdots & \ddots & \vdots \\
       x_{k+1}^{n-k-1}+P_1x_{k+1}^{n-k-2}+\cdots+P_{n-k-1} & \cdots &
       x_{n+1}^{n-k-1}+P_1x_{n+1}^{n-k-2}+\cdots+P_{n-k-1} \\
       \big[x_1,\dots,x_k,x_{k+1};f \big]& \cdots
	  & \big[x_1,\dots,x_k,x_{n+1};f \big]
    \end{array}\right| \\[2mm]
    &= \left|\begin{array}{ccc}
    1 & \cdots & 1 \\
    x_{k+1} & \cdots & x_{n+1} \\
    \vdots & \ddots & \vdots \\
    x_{k+1}^{n-k-1} & \cdots & x_{n+1}^{n-k-1} \\
    \big[x_1,\dots,x_k,x_{k+1};f \big]& \cdots
	  & \big[x_1,\dots,x_k,x_{n+1};f \big]
    \end{array}\right|.
}
Replacing the right hand side of \eq{fid} by the right hand side the above identity, it follows from 
\thm{2} that the convexity of $f$ with respect to the Chebyshev system 
$(\omega_1,\dots,\omega_n)=(p_0,\dots,p_{n-1})$ is equivalent to the monotonicity/convexity properties of 
the mapping $x\mapsto \big[ x_1,\dots,x_k,x;f \big]$ on $ H\setminus\{ x_1,\dots,x_k \} $ or on the 
subintervals of $ H\setminus\{ x_1,\dots,x_k \} $.
\end{proof}

\Exa{*}{For a direct application of our results, let $ \omega$ be the $3$-dimensional positive Chebyshev 
system $\big(1,\cos,\sin \big) $ over the interval $H=]-\pi,0[$. Then $\omega_{\<2\>}=(1,\cos)$ is a 
$2$-dimensional positive Chebyshev system and, for any function $f:H\to\R $, the following statements are 
equivalent
\begin{enumerate}[(i)]
 \item $f$ is $\omega$-convex on $H$; 
 \item For each $x_1\in H$, the function $ x\mapsto \big[x_1,x;f \big]_{(1,\cos)} = 
\frac{f(x)-f(x_1)}{\cos(x)-\cos(x_1)}$ is convex on $H\setminus\{x_1\}$ with respect to Chebyshev system $ 
\Big(1,-\ctg\big(\frac{x_1+(\cdot)}{2}\big) \Big) $.
\end{enumerate}}

\begin{proof}
By well-known trigonometrical identities
\Eq{*}{
\sin(x)-\sin(y)=2\cos\Big(\frac{x+y}{2}\Big) \sin\Big(\frac{x-y}{2}\Big),\qquad
\cos(x)-\cos(y)=-2\sin\Big(\frac{x+y}{2}\Big) \sin\Big(\frac{x-y}{2}\Big),
}
for $x,y\in H$ with $x\neq y$, we have
\Eq{id+}{
  \frac{\sin(x)-\sin(y)}{\cos(x)-\cos(y)}=-\ctg\Big(\frac{x+y}{2}\Big).
}
Using \eq{id+}, the right hand side of \eq{fid} with $k=1$ can be written as
\Eq{*}{
      &\left|\begin{array}{ccc}
    1 & 1 & 1 \\
    \big[x_1,x_2;\sin \big]_{(1,\cos)} & 
       \big[x_1,x_3;\sin \big]_{(1,\cos)} & 
       \big[x_1,x_4;\sin \big]_{(1,\cos)} \\
    \big[x_1,x_2;f \big]_{(1,\cos)}&  \big[x_1,x_3;f \big]_{(1,\cos)}
        & \big[x_1,x_4;f \big]_{(1,\cos)}
    \end{array}\right|\\
    &\qquad\qquad=\left|\begin{array}{ccc}
    1 & 1 & 1 \\
    \frac{\sin(x_2)-\sin(x_1)}{\cos(x_2)-\cos(x_1)}& 
       \frac{\sin(x_3)-\sin(x_1)}{\cos(x_3)-\cos(x_1)} & 
       \frac{\sin(x_4)-\sin(x_1)}{\cos(x_4)-\cos(x_1)} \\
    \big[x_1,x_2;f \big]_{(1,\cos)}&  \big[x_1,x_3;f \big]_{(1,\cos)}
        & \big[x_1,x_4;f \big]_{(1,\cos)}
    \end{array}\right|\\
    &\qquad\qquad=\left|\begin{array}{ccc}
    1 & 1 & 1 \\
    -\ctg{\big(\frac{x_1+x_2}{2}\big)} & 
       -\ctg{\big(\frac{x_1+x_3}{2}\big)} & 
       -\ctg{\big(\frac{x_1+x_4}{2}\big)} \\
    \big[x_1,x_2;f \big]_{(1,\cos)}&  \big[x_1,x_3;f \big]_{(1,\cos)}
        & \big[x_1,x_4;f \big]_{(1,\cos)}
    \end{array}\right|.
}
By \thm{2}, the $(1,\cos,\sin)$-convexity of a function $f:H\to\R$ is equivalent to the nonnegativity of the 
above determinants, which exactly means that the function $ x\mapsto\big[x_1,x;f \big]_{(1,\cos)} $ is convex 
with respect to Chebyshev system $ \Big(1,\ctg\big(\frac{-x_1-(\cdot)}{2}\big) \Big) $.
\end{proof}

\section{Differences of $\omega$-convex functions}

Let $H$ be an open real interval throughout this section and let $\omega:H\to\R^n$ be a Chebyshev system.
We introduce the notion of $\omega$-variation which will turn out to be finite for differences of 
$\omega$-convex functions.

Given a subinterval $[a,b]$, define the set of partitions $\P([a,b])$ of $[a,b]$ by 
\Eq{*}{
  \P([a,b]):=\{(x_0,\dots,x_n)\mid n\in\N,\,a=x_0<\cdots<x_n=b\}.
}
The \emph{$\omega$-variation} of $f:H\to\R$ on $[a,b]$ is now defined by
\Eq{*}{
  V^\omega_{[a,b]}(f)
  :=\sup\bigg\{\sum_{i=0}^{m-n}\Big|\big[x_{i+1},\dots,x_{i+n};f \big]_{\omega}
       -\big[x_i,\dots,x_{i+n-1};f \big]_{\omega}\Big|
       \,:\,m\geq n,\,(x_0,\dots,x_m)\in\P([a,b])\bigg\}.
}

Applying the notion of $\omega$-variation, the next theorem formulates a necessary condition in order that a 
function could be decomposed as the difference of two $\omega$-convex functions. The sufficiency of this 
conditions remains an open problem.

\Thm{3}{Let $\omega:=(\omega_1,\dots,\omega_n):H\to \R^n$ be an $n$-dimensional positive Chebyshev system 
such that $\omega_{\<n-1\>}$ is an $(n-1)$-dimensional positive Chebyshev system and let $f:H\to \R$ be a 
function. If there exist $\omega$-convex functions $g,h:H\to\R$ such that $f=g-h$, then, for all subinterval 
$[a,b]\subseteq H$, the $\omega$-variation $V^\omega_{[a,b]}(f)$ is finite. Furthermore, for all elements 
$a_1,\dots,a_n,b_1,\dots,b_n\in H$ with $a_1<\dots<a_n=a$ and $b=b_1<\dots<b_n$, the inequality
\Eq{V}{
  V^\omega_{[a,b]}(f)\leq\big[b_1,\dots,b_n;g+h \big]_{\omega}-\big[a_1,\dots,a_n;g+h \big]_{\omega}
}
holds.}

\begin{proof} Assume that $f$ is of the form $f=g-h$, where $g,h:H\to\R$ are $\omega$-convex functions.
Let $a,b\in H$ with $a<b$ and fix $a_1<\dots<a_n=a$ and $b=b_1<\dots<b_n$ in $H$.
Let $(x_0,\dots,x_m)\in\P([a,b])$ be an arbitrary partition of $[a,b]$ with $m\geq n$.
Then, using the linearity of $\omega$-divided differences and the triangle inequality, we get  
\Eq{*}{
    &\sum_{i=0}^{m-n}\Big|\big[x_{i+1},\dots,x_{i+n};f \big]_{\omega}
         -\big[x_i,\dots,x_{i+n-1};f \big]_{\omega}\Big|\\&=  
    \sum_{i=0}^{m-n}\Big|\big[x_{i+1},\dots,x_{i+n};g-h \big]_{\omega}
           -\big[x_i,\dots,x_{i+n-1};g-h \big]_{\omega}\Big| \\ &\leq
    \sum_{i=0}^{m-n}\bigg(\Big|\big[x_{i+1},\dots,x_{i+n};g \big]_{\omega}
	-\big[x_i,\dots,x_{i+n-1};g \big]_{\omega}\Big|+
	\Big|\big[x_{i+1},\dots,x_{i+n};h \big]_{\omega}
	-\big[x_i,\dots,x_{i+n-1};h \big]_{\omega}\Big|\bigg).
}
In view of the monotonicity property of $\omega$-divided differences established in \cor{1} for 
the $\omega$-convex functions $g$ and $h$, for $i\in\{0,\dots,m-n\}$, we have
\Eq{*}{
  \Big|\big[x_{i+1},\dots,x_{i+n};g \big]_{\omega}
	-\big[x_i,\dots,x_{i+n-1};g \big]_{\omega}\Big|
  &=\big[x_{i+1},\dots,x_{i+n};g \big]_{\omega}
	-\big[x_i,\dots,x_{i+n-1};g \big]_{\omega}, \\
  \Big|\big[x_{i+1},\dots,x_{i+n};h \big]_{\omega}
	-\big[x_i,\dots,x_{i+n-1};h \big]_{\omega}\Big|
  &=\big[x_{i+1},\dots,x_{i+n};h \big]_{\omega}
	-\big[x_i,\dots,x_{i+n-1};h \big]_{\omega}.
}
Thus, performing telescopic summation and using the linearity of $\omega$-divided differences, we obtain
\Eq{*}{
    &\sum_{i=0}^{m-n}\bigg(\Big|\big[x_{i+1},\dots,x_{i+n};g \big]_{\omega}
	-\big[x_i,\dots,x_{i+n-1};g \big]_{\omega}\Big|+
	\Big|\big[x_{i+1},\dots,x_{i+n};h \big]_{\omega}
	-\big[x_i,\dots,x_{i+n-1};h \big]_{\omega}\Big|\bigg)\\
    &=\sum_{i=0}^{m-n}\bigg(\big[x_{i+1},\dots,x_{i+n};g \big]_{\omega}
	-\big[x_i,\dots,x_{i+n-1};g \big]_{\omega}	
	+\big[x_{i+1},\dots,x_{i+n};h \big]_{\omega}
	-\big[x_i,\dots,x_{i+n-1};h \big]_{\omega}\bigg)\\
    &=\bigg(\big[x_{m-n+1},\dots,x_m;g \big]_{\omega}-
         \big[x_0,\dots,x_{n-1};g \big]_{\omega}\bigg)
     +\bigg(\big[x_{m-n+1},\dots,x_m;h \big]_{\omega}-
         \big[x_0,\dots,x_{n-1};h \big]_{\omega}\bigg)\\
   &=\big[x_{m-n+1},\dots,x_m;g+h \big]_{\omega}
         -\big[x_0,\dots,x_{n-1};g+h \big]_{\omega}.
}
Now, using the inequalities $a_i<x_{i-1}$ and $x_{m-n+i}<b_i$ (which follow from the choice of 
$a_1,\dots,a_n$ and $b_1,\dots,b_n$), and applying again the monotonicity property of $\omega$-divided 
differences established in \cor{1} for the $\omega$-convex function $g+h$, we get
\Eq{*}{
   -\big[x_0,\dots,x_{n-1};g+h \big]_{\omega}&\leq -\big[a_1,\dots,a_n;g+h \big]_{\omega}\\
    \big[x_{m-n+1},\dots,x_m;g+h \big]_{\omega}&\leq \big[b_1,\dots,b_n;g+h \big]_{\omega}.
}
Finally, combining the above three estimates, for every partition $(t_0,\dots,t_m)$ of $[a,b]$, we obtain
\Eq{*}{
  \sum_{i=0}^{m-n}\Big|\big[x_{i+1},\dots,x_{i+n};f \big]_{\omega}
         -\big[x_i,\dots,x_{i+n-1};f \big]_{\omega}\Big|
  \leq \big[b_1,\dots,b_n;g+h \big]_{\omega}-\big[a_1,\dots,a_n;g+h \big]_{\omega},
}
which implies
\Eq{*}{
  V^\omega_{[a,b]}(f)\leq\big[b_1,\dots,b_n;g+h \big]_{\omega}-\big[a_1,\dots,a_n;g+h \big]_{\omega} <+\infty.
}
Thus, inequality \eq{V} and the theorem is proved.
\end{proof}

For the case of higher-order convexity in the sense of Hopf and Popoviciu, the following characterization 
holds (cf.\ \cite{Pop44}), which, in one direction, is a consequence \thm{3}. 

\Cor{3}{
Let $\omega:=(p_0,\dots,p_{n-1}):H\to\R^n$ be the $n$-dimensional positive Chebyshev system and let $f:H\to \R$ be a 
function. Then, there exist $\omega$-convex functions $g,h:H\to\R$ such that $f=g-h$ if and only if for all subinterval 
$[a,b]\subseteq H$, the $\omega$-variation $V^\omega_{[a,b]}(f)$ is finite.}

%\bibliography{publ,funcequ}

\begin{thebibliography}{1}

\bibitem{AkrAkrMal96}
A.~G. Akritas, E.~K. Akritas, and G.~I. Malaschonok, \emph{{Various proofs of
  {S}ylvester's (determinant) identity}}, Math. Comput. Simulation \textbf{42}
  (1996), no.~4-6, 585--593, Symbolic computation, new trends and developments
  (Lille, 1993). \MR{1430843}

\bibitem{BesPal04}
M.~Bessenyei and Zs. P\'ales, \emph{{On generalized higher-order convexity and
  {H}ermite--{H}adamard-type inequalities}}, Acta Sci. Math. (Szeged)
  \textbf{70} (2004), no.~1-2, 13--24. \MR{2005e:26012}

\bibitem{Gan59}
F.~R. Gantmacher, \emph{{The theory of matrices. {V}ols. 1, 2}}, {Translated by
  K. A. Hirsch}, Chelsea Publishing Co., New York, 1959. \MR{0107649}

\bibitem{Hop26}
E.~Hopf, \emph{{{\"U}ber die {Z}usammenh\"ange zwischen gewissen h\"oheren
  {D}ifferenzenquotienten reeller {F}unktionen einer reellen {V}ariablen und
  deren {D}ifferenzierbarkeitseigenschaften}}, Ph.D. thesis,
  Friedrich--Wilhelms--Universit\"at Berlin, 1926.

\bibitem{Kar68}
S.~Karlin, \emph{{Total positivity. {V}ol. {I}}}, Stanford University Press,
  Stanford, California, 1968. \MR{37 \#5667}

\bibitem{KarStu66}
S.~Karlin and W.~J. Studden, \emph{{Tchebycheff systems: {W}ith applications in
  analysis and statistics}}, {Pure and Applied Mathematics, Vol. XV},
  Interscience Publishers John Wiley \& Sons, New York-London-Sydney, 1966.
  \MR{34 \#4757}

\bibitem{Pop44}
T.~Popoviciu, \emph{{Les fonctions convexes}}, Hermann et Cie, Paris, 1944.
  \MR{8,319a}

\bibitem{PalRad16}
Zs. P\'ales and \'E.~Sz. Rad\'acsi, \emph{{Characterizations of higher-order
  convexity properties with respect to {C}hebyshev systems}}, Aequationes Math.
  \textbf{90} (2016), no.~1, 193--210. \MR{3471295}

\bibitem{Was06b}
Sz. W\c{a}sowicz, \emph{{Some properties of generalized higher-order convexity}},
  Publ. Math. Debrecen \textbf{68} (2006), no.~1-2, 171--182. \MR{2213549 (2006m:26016)}

\end{thebibliography}
%\bibliographystyle{amsplain}

\def\MR#1{}

\end{document}